\documentclass{amsart}
\usepackage{amsmath,amssymb,bbm,epsfig}

\newcommand{\Isom}[1]{\mathop{Isom}(#1)}
\newcommand{\Tree}{\mathcal{T}}

\newcommand{\Sym}[1]{S_{#1}}

\theoremstyle{plain}
\newtheorem{theorem}{Theorem}
\newtheorem{proposition}{Proposition}
\newtheorem{lemma}{Lemma}
\newtheorem{corollary}{Corollary}

\theoremstyle{definition}
\newtheorem{definition}{Definition}
\newtheorem*{remark}{Remark}
\newtheorem*{notation}{Notation}

\begin{document}

\title{Homomorphic images of Branch groups, and Serre's property (FA).}

\author{Thomas Delzant}
\address{Thomas Delzant, D\'epartement de math\'ematiques, 7 rue Descartes, Universit\'e de Strasbourg, F-67084
Strasbourg.} \email{delzant@math.u-strasbg.fr}

\author{Rostislav Grigorchuk}
\address{Rostislav Grigorchuk, Department of Mathematics, Texas A \&M University, MS-3368,
College Station, TX, 77843-3368, USA} \email{grigorch@math.tamu.edu}
\thanks{Partially supported by NSF grants DMS-0600975 and DMS-0456185}

\begin{abstract}
It is shown that a finitely generated branch group has Serre's
property (FA) if and only if it does not surject onto the infinite
cyclic group or the infinite dihedral group. An example of a
finitely generated self-similar branch group surjecting onto the
infinite cyclic group is constructed.
\end{abstract}

\dedicatory{Dedicated to the memory of Sasha Reznikov}

\date{May 2, 2007}

\maketitle

\section*{Introduction}

The study of groups acting on trees is a central subject in
geometric group theory. The Bass-Serre theory establishes a
dictionary between the geometric study of groups acting on trees and
the algebraic study of amalgams and HNN extensions. A central topic
of investigation is the fixed point property for groups acting on
trees, introduced by J.-P. Serre in his book as the property
(FA)\cite{serre:trees}. A fundamental result due to Tits states that
a group without a free subgroup on two generators which acts on a
tree by automorphisms fixes either a vertex or a point on the
boundary or permutes a pair of points on the boundary; see
\cite{tits:lie_kolchin, PaysValette}. The group $SL(3,\mathbb{Z})$,
and more generally, groups with Kazdhan's property (T), in
particular lattices in higher rank Lie groups have the property (FA)
(\cite{harpe_valette, margulis:lie}). A natural problem is to
understand the structure of the class of (FA)-groups (the class of
groups having the property (FA)). There is an algebraic
characterization of enumerable (FA)-groups, due to J.-P. Serre.
(\cite{serre:trees}, Theorem I.6.15, page 81).

An enumerable group has the property (FA) if and only if it
satisfies the following three conditions:

(i) it is not an amalgam,

 (ii) it is not
indicable (i.e. admits no epimorphism onto $\mathbb{Z}$),

(iii) it is finitely generated.

But even such a nice result does not clarify the structure of the
class of (FA)-groups, as the first of these properties is usually
difficult to check.

The class of (FA)-groups contains the class of finite groups and is
closed under quotients. As every infinite finitely generated group
surjects onto a just-infinite group (i.e. an infinite group with all
proper quotients finite) a natural problem is to describe just
infinite (FA)-groups.

In~\cite{gri:branch} the class (JINF) of just infinite groups is
divided in three subclasses: the class (B) of branch groups, the
class (HJINF) of finite extensions of finite powers of hereditary
just infinite groups and the class (S) of finite extensions of
finite powers of simple groups. For example, the group
$SL(3,\mathbb{Z})$ belongs to the class (JINF); all infinite
finitely generated simple torsion groups constructed
in~\cite{olshansk:inf_noeth} are (FA)-groups and belong to the class
(S).

A precise definition of a branch group is given in Section 1.
Roughly speaking a branch group is a group which acts faithfully and
level transitively on a spherically homogeneous rooted tree, and for
which the structure of the lattice of subnormal subgroups mimics the
structure of the tree. Branch groups may enjoy unusual properties.
Among them one can find finitely generated infinite torsion groups,
groups of intermediate growth, amenable but not elementary amenable
groups and other surprising objects. Profinite branch groups are
also related to Galois theory and other topics in Number Theory
\cite{boston:00}.

In this article we discuss fixed point properties for actions of
branch groups on Gromov hyperbolic spaces, in particular on
$\mathbb{R}$-trees. One of the corollaries of Theorem 3 is:

\begin{theorem}\label{A}
Let $G$ be a finitely generated branch group. Then $G$ has Serre's
property (FA) if and only if it is not indicable and has no
epimorphism onto $\mathbbm D_\infty$.
\end{theorem}

Therefore a branch group cannot be an amalgam unless it surjects
onto $\mathbbm D_\infty$.

We say that a group is (FL) if it has no epimorphism onto
$\mathbb{Z}$ or $\mathbbm{D}_{\infty}$. A f.g. group is (FL) if and
only if it fixes a point whenever it acts isometrically on a line.

All proper quotients of branch groups are virtually
abelian~\cite{gri:branch}. A quotient of a branch group may be
infinite: the full automorphism group of the binary rooted tree is a
branch group and its abelianization is the infinite cartesian
product of copies of a group of order two. It is more difficult to
construct examples of finitely generated branch groups with infinite
quotients (especially in the restricted setting of self-similar
groups). The corresponding question was open since 1997 when the
second  author introduced the notion of a branch group. Perhaps the
main difficulty was psychological, as he (and some other researches
working in the area) was sure that all finitely generated branch
groups are just infinite. Now we know that this is not correct and
the second part of the paper (Section 3) is devoted to a
construction of an example of an indicable finitely generated branch
group (thus providing an example of a finitely generated branch
group without the property (FA)). This example is the first example
of a finitely generated branch group defined by a finite automaton
that is not just infinite. Another example is related to Hanoi
Towers group on 3 pegs $H$ (introduced
in~\cite{grigorchuk-s:hanoi-cr} and independently
in~\cite{nekrash:book}). Hanoi Towers group $H$ is a 3-generated
branch group~\cite{grigorchuk-s:standrews} that has a subgroup of
index 4 (the Apollonian group) which is also a branch group and is
indicable (this is announced in\cite{gri_nek_sun:hanoi}). The group
$H$ itself is not indicable (it has finite abelianization), but it
surjects onto $\mathbbm D_\infty$, as was recently observed by Zoran
\v{S}uni\'c. Thus $H$ is the first example of a finitely generated
branch group defined by a finite automaton that surjects onto
$\mathbbm D_\infty$.

The example presented in this paper is an elaboration of the
3-generated torsion 2-group $G = \langle a,b,c,d \rangle$ firstly
constructed in~\cite{gri:80_en} and later studied in~\cite{gri:hilb,
grig_maki:93,gri:example,gri:schur} and other papers (see also the
Chapter VIII of the book~\cite{harpe} and the article
\cite{CecchMachiScar}.

Let $L$ be the group generated by the automaton defined in
Figure~\ref{D1}.

\begin{theorem}\label{B}
The group $L$ is a branch, contracting group that surjects onto
$\mathbb{Z}$.
\end{theorem}

An interesting question is to understand which virtually abelian
group can be realized as a quotient of a finitely generated branch
group. This question is closely related to the problem of
characterization of finitely generated branch groups having the
Furstenberg-Tychonoff fixed ray property (FT)~\cite{grig:tychon})
(existence of an invariant ray for actions on a convex cone with
compact base).

\subsection*{Acknowledgments}
The authors are thankful to Zoran \v{S}uni\'c and Laurent Bartholdi
for valuable discussions, comments, and suggestions.


\section{Basic definitions and some notations.}
Let $T$ be a tree, $G$ be a group acting on $T$ (without inversion
of edges) and $T^G$ the the set of fixed vertices of $T$.

\begin{definition} A group $G$ has the property (FA) if
for every simplicial tree $T$ on which $G$ acts simplicially and
without inversion, $T^G\neq\emptyset$.
\end{definition}

The class of (FA)-groups possesses the following properties.

(i) The class of (FA)-groups is closed under taking quotients.

(ii) Let $G$ be a group with the property (FA). If $G$ is a subgroup
of an amalgamated free product $G_1*_{A}G_2$ or an HNN extension
$G=G_1*_A$, then $G$ is contained in a conjugate of $G_1$ or $G_2$.

(iii) The class of (FA)-groups is closed under forming extensions.

(iv) If a subgroup of finite index in a group $G$ has the property
(FA), then the group $G$ itself has the property (FA).

(v) Every finitely generated torsion group has the property (FA).

The class of (FA)-groups has certain nice structural properties and
is interesting because of the strong embedding property given by
(ii) and by the fact that the eigenvalues of matrices in the image
of a linear representation $\rho:G \rightarrow GL_2(\textit{k})$ are
integral over $\mathbb{Z}$ for any field $\textit{k}$ (Prop. 22,
\cite{serre:trees}).

The property (i), the existence of just infinite quotients for
finitely generated infinite groups and the trichotomy from
\cite{gri:branch} mentioned in the introduction make the problem of
classification of finitely generated just infinite (FA)-groups
worthwhile. We are reduced to the classification of finitely
generated (FA)-groups in each of the classes (B), (HJINF) and (S).
Below we solve this problem, in a certain sense, for the class (B).

If a group $G$ has a quotient isomorphic to $\mathbb{Z}$, then it
acts by translations on a line and cannot be an (FA)-group.
Similarly, if $G$ surjects onto the infinite dihedral group
$\mathbbm D_{\infty}$, then it acts on the line via the obvious
action of $\mathbbm D_{\infty}$. This suggests the following
definition (the first part being folklore):

\begin{definition}
a)A group is called \emph{indicable} if it admits an epimorphism
onto $\mathbb{Z}$.

b) A group has property (FL) (fixed point on line) if every action
of $G$ by isometry on a line fixes a point. If $G$ is finitely
generated this means that $G$ has no epimorphism onto $\mathbb Z$ or
$\mathbbm D _\infty$.
\end{definition}

In this article we will often use two other notions: the notion of a
hyperbolic space and that of a branch group.

For the definition and the basic properties of Gromov hyperbolic
spaces we refer the reader to~\cite{coo_del_pap}. The theory of
CAT$(0)$-spaces is described in~\cite{bri_hae}. For the definition
and the study of basic properties of branch groups we refer the
reader to~\cite{gri:branch,bar_gri_sun:branch}.

Let us recall the main definition and a few important facts and
notations that will be often used later.

\begin{definition}
 A group $\mathcal{G}$ is an
\emph{algebraically branch group} if there exists a sequence of
integers $\overline{k}=\{k_n\}_{n=0}^\infty$ and two decreasing
sequences of subgroups $\{R_n\}_{n=0}^\infty$ and
$\{V_n\}_{n=0}^\infty$ of $\mathcal{G}$ such that

(1) $k_n \geq 2$, for all $n>0$, $k_0=1$,

(2) for all $n$,
\begin{equation}\label{branchdecomp}
 R_n = V_n^{(1)} \times V_n^{(2)} \times \dots \times V_n^{(k_0k_1 \dots k_n)},
\end{equation}
where each $V_n^{(j)}$ is an isomorphic copy of $V_n$,

(3) for all $n$, the product decomposition \eqref{branchdecomp} of
$R_{n+1}$ is a refinement of the corresponding decomposition of
$R_n$ in the sense that the $j$-th factor $V_n^{(j)}$ of $R_n$,
$j=1,\dots,k_0k_1\dots k_n$ contains the $j$-th block of $k_{n+1}$
consecutive factors
\[ V_{n+1}^{((j-1)k_{n+1}+1)} \times \dots \times V_{n+1}^{(jk_{n +1})}\]
of $R_{n+1}$,

(4) for all $n$, the groups $R_n$ are normal in $\mathcal{G}$ and
\[ \bigcap_{n=0}^\infty R_n = 1, \]

(5) for all $n$, the conjugation action of $\mathcal{G}$ on $R_n$
permutes transitively the factors in~\eqref{branchdecomp},

and

(6) for all $n$, the index $[\mathcal{G}:R_n]$ is finite.

A group $G$ is a \emph{weakly algebraically branch group} if there
exists a sequence of integers $\overline{k}=\{k_n\}_{n=0}^\infty$
and two decreasing sequences of subgroups $\{R_n\}_{n=0}^\infty$ and
$\{V_n\}_{n=0}^\infty$ of $G$ satisfying the conditions (1)-(5).
\end{definition}

There is a geometric counterpart of this definition.

Let $(\mathcal{T}, \emptyset)$ be a spherically homogeneous rooted
tree, where $\emptyset $ is the root and $\mathcal{G}$ be a group
acting on $(\mathcal{T}, \emptyset)$ by automorphisms preserving the
root. Let $v$ be a vertex, and $\mathcal{T}_v$ be the subtree
consisting of the vertices $w$ such that $v\in [w,\emptyset]$
(geodesic segment joining $w$ with the root). The rigid stabilizer
$rist_\mathcal{G}(v)$ of a vertex $v$ consists of elements acting
trivially on $\Tree \setminus \Tree_v$. The rigid stabilizer of the
$n$-th level, denoted $rist_\mathcal{G}(n)$, is the group generated
by the rigid stabilizers of the vertices on level $n$.

The action of $\mathcal G$ on $\mathcal T$ is called geometrically
branch if it is faithfull, level transitive, and if, for any $n$,
the rigid stabilizer $rist_\mathcal{G}(n)$ of $n-$th level of the
tree has finite index in $\mathcal{G}$.

Observe that, in the level transitive case, the rigid stabilizers of
the vertices of the same level are conjugate in $\mathcal{G}$. In
this case $rist_\mathcal{G}(n)$ is algebraically isomorphic to the
product of copies of the same group (namely the rigid stabilizer of
any vertex on the given level). Hence the rigid stabilizers of the
levels and vertices play the role of the subgroups $R_n$ and $V_n$
of the algebraic definition. A geometrically branch group is
therefore algebraically branch. The algebraic definition is slightly
more general than the geometric one but at the moment it is not
completely clear how big the difference between the two classes of
groups is. Observe that in Section 2 we will assume that the
considered groups are algebraically branch, while in sections 3 and
4 we construct examples of geometrically branch groups.

When constructing these examples, we will deal only with actions on
a rooted binary tree and our notation and the definition below are
adapted exactly for this case. Let $\mathcal{G}$ be a branch group
acting on a binary rooted tree $\mathcal{T}$. The vertices of
$\mathcal{T}$ are labeled by finite sequences of $0$ and $1$. Let
$\mathcal{T}_0,\mathcal{T}_1$ be the two subtrees consisting of the
vertices starting with $0$ or $1$, respectively.

\begin{notation}\label{2.6}
If $A, B, C \subset \mathop{Aut} (\mathcal{T})$ are three subgroups,
we write $A \succeq B \times C$ if $A$ contains the subgroup $B
\times C$ of the product \text{$\mathop{Aut} (\mathcal{T}_0) \times
\mathop{Aut} (\mathcal{T}_1)$ }via the canonical identification of
$\mathop{Aut} (\mathcal{T})$ with $\mathop{Aut} (\mathcal{T}_i)$.
\end{notation}

Recall that a level transitive group $\mathcal{G}$ acting on a
regular rooted binary tree is called regular branch over its normal
subgroup $H$ if $H$ has finite index in $\mathcal{G}$, $H \succeq H
\times H$ and if moreover the last inclusion is of finite index.

A level transitive group $\mathcal{G}$ is called weakly regular
branch over a subgroup $H$ if $H$ is nontrivial and $H \succeq H
\times H$.

\begin{definition}\label{2.4}
A group $\mathcal{G}$ acting on the rooted binary tree
$(\mathcal{T}, \emptyset)$ is called \emph{self-replicating} if, for
every vertex $u$, the image of the stabilizer $s t_\mathcal{G} ( u
)$ of $u$ in $\mathop{Aut}(\Tree_u)$ (the automorphism group of the
rooted tree $\mathcal{T}_u)$ coincides with the group $G$ after the
canonical identification of $\Tree$ with $\Tree_u$.
\end{definition}

Obviously a self-replicating group is level transitive if and only
if it is transitive on the first level (see also Lemma~A in
\cite{gri:branch}).

We will use the notations $\langle R \rangle^\mathcal{G}$ for the
normal closure in $\mathcal{G}$ of a subset $R
\subset{\mathcal{G}}$, $x^y=y^{-1}xy$, $[x,y]=x^{-1}y^{-1}xy$. Given
two subgroups, $A,B$ in a group $\mathcal G$, $[A,B]$ is the
subgroup of $\mathcal G$ generated by the commutators $[a,b]$ of
elements in $A$ and $B$, and $[A,B,]^{\mathcal G}$ its normal
closure. If $\mathcal{G}$ is a group, $\gamma_2( \mathcal{G})$
denote the seond member of its lower central series.

\section{Fixed point properties of branch groups}

Let $X$ be a Gromov hyperbolic metric space, and $\partial X$ its
Gromov's boundary. Recall (see~\cite{Gromov:Hyperbolic} or
\cite{coo_del_pap} chap. 9 for instance) that a subgroup $G$ of the
group $\Isom{X}$ of isometries of $X$ is called elliptic if it has a
bounded orbit (or equivalently if every orbit is bounded), parabolic
if it has a unique fixed point on $\partial X$ but is not elliptic,
and loxodromic if it is not elliptic and if there exists a pair
$w^+, w^-$ of points in $\partial X$ preserved by $G$. A group which
is either elliptic, or parabolic or loxodromic is called elementary;
this terminology is inspired by the theory of Kleinian groups. There
are no constraints on the algebraic structure of elementary groups
due to the following remark.

\begin{remark}\label{1.2}
Every f.g. group $G$ can be realized as a \textit{parabolic} group
of isometries of some proper geodesic hyperbolic space: if $C$ is
the Cayley graph of $G$, $C \times \mathbbm{R}$ admits a
$G$-invariant hyperbolic metric (~\cite{Gromov:Hyperbolic}, 1.8.A,
note that this construction is equivariant). One can also construct
a finitely generated group acting on a tree with a unique fixed
point at infinity. For instance the lamplighter group (semi-direct
product of $\mathbbm Z $ and $\mathbbm Z_2[t,t^{-1}]$) fixes a
unique point in the boundary of the tree of $GL_2( \mathbbm
Z_2[t^{-1},t])$. In fact, the lamplighter group can be indentified
with upper triangular matrices with one eigenvalue equal to $1$ the
other being $t^n$. As all these matrices have a common eigenspace,
they fix one point in the boundary of the tree of $GL_2$ (the
projective line on $\mathbbm Z_2[t,t^{-1}]]$); but this group
contains the Jordan matrix and therefore cannot fix two points in
the boundary of this tree.

\end{remark}

In what follows, $X$ denotes a complete Gromov hyperbolic geodesic
space. We will assume that either $X$ is proper (closed balls are
compact) or that $X$ is a complete $\mathbbm{R}$-tree, i.e. a
complete $0$-hyperbolic geodesic metric space. In the first case, $X
\cup \partial X$ is a compact set (in the natural topology) and an
unbounded sequence of points in $X$ admits a subsequence which
converges to a point in $\partial X$. Important examples of such
spaces are Cayley graphs of hyperbolic groups (see
\cite{Gromov:Hyperbolic} for instance). Other examples are universal
covers of compact manifolds of non positive curvature. Note that
properness implies completness for a metric space, but the converse
is false. Recall also that a geodesic space is proper if and only if
it is complete and locally compact~\cite{gromov:metric_struct}. The
Gromov hyperbolicity of a geodesic space can be defined in several
ways (thineness of geodesic triangles, properties of the Gromov
product etc.) which are equivalent (see~\cite{coo_del_pap} Chap.~1);
we will prefer the definition in terms of the Gromov product (
\cite{coo_del_pap} chap.1, def. 1.1).

For the rest of the statements in this section we will assume
 that the following condition on the pair $(X,G)$ holds:

 (C) $X$ is a complete geodesic space and $X$ is either
proper hyperbolic or is an $\mathbbm{R}$-tree. $G$ is a group and
$\varphi: G \rightarrow \Isom{X}$ is an isometric action of $G$ on
$X$.

Note that such an action extends uniquely to a continuous action on $X
\cup
\partial X$.

\begin{theorem}\label{1.1} Let $G$ be a branch group acting isometrically on a hyperbolic space $X$. Suppose the pair $(X,G)$ satisfies the
condition (C).
 Then

a) the image of $G$ in $\Isom{X}$ is elementary.

b) Suppose furthermore that G satisfies the property (FL), and $X$
is a hyperbolic graph with uniformly bounded valence of vertices.
Then $\varphi(G)$ is elliptic or parabolic.

c) If $X$ is CAT(0) and if the group $\varphi(G)$ is elliptic then
it has a fixed point in $X$.

 d) If X is CAT(-1), or is an
$\mathbbm{R} -$tree then $\varphi(G)$ fixes a point
 in $X$ or in $\partial X$, or preserves a line in $\ X$.

 e) Let X be an
$\mathbbm{R} -$tree. Suppose further that $G$ is f.g.; then $G$
cannot be parabolic.
\end{theorem}

\begin{corollary}\label{cor:1}
Let $G$ be a f.g. branch group. $G$ has fixed point property for
actions on $\mathbbm R$-trees if and only if it has property (FL).
\end{corollary}
\begin{proof}
A tree is $CAT(-1)$, so if $G$ acts on a tree and does not fix a
point, it must either preserve a line or a unique point on $\partial
X$. The last possibility is excluded by $e)$.
\end{proof}

>From $d)$ we also deduce:
\begin{corollary}\label{cor:2}
If $X$ is $CAT(-1)$ and $G$ acts on $X \cup \partial X$ by
isometries, then $G$ fixes a point or contains a subgroup of index 2
which fixes two points in $\partial X$.
\end{corollary}

Recall (see~\cite{serre:trees}) that a group $G$ is an amalgam
(resp. an HNN extension) if it can be written as a free product with
amalgamation $G=A*_CB$, with $C\not =A,B$ (resp.
$G=A*_{tCt^-1=C'})$. We say that this amalgam (resp. HNN extension)
is $strict$ if the index of $C$ in $A$ is at least 3 and the index
of $C$ in $B$ is at least $2$ (resp. the indexes of $C$ and $C'$ in
$A$ are at least 2). If $G$ splits as an amalgam or HNN extension
then $G$ acts on a simplicial tree $T$ without edge inversion s.t.
$T/G$ has one edge and 2 vertices in the case of an amalgam, and one
edge and one vertex in the case of an HNN extension. It is easy too
see that if a group is a strict amalgam or HNN extension its action
on Serre's tree is not elementary. If $G=A*_CB$ with $C$ of index 2
in $A$ and $B$ Serre's tree is a line, and $G$ permutes the two ends
of this line. If $G=A*_{tCt^-1=C'}$ and $C=C'=A$, Serre's tree is a
line and $G$ fixes the two ends of this line. If $G=A*_{tCt^-1=C'}$
is a strictly ascending HNN extension ($C'=A$, but $C \not =A$), the
group $G$ contains a hyperbolic element (the letter $t$ for
instance) and fixes exactly one end of the the tree. Therefore the
property $e)$ implies the following:

\begin{corollary}\label{cor:3}
Let $G$ be a f.g. branch group. Then $G$ is neither a strict amalgam
nor a strict HNN extension nor a strictly ascending HNN
extension.
\end{corollary}

Before proving Theorem~\ref{1.1} let us state and prove some
statements that have independent interest and will be used later.

 Recall that an isometry $f$ of a hyperbolic space $X$ is called
elliptic (resp. parabolic, resp. hyperbolic) if the subgroup
generated by$ f$ is elliptic (resp. parabolic, resp. loxodromic). It
can be proved (see~\cite{coo_del_pap}, chap.~9) that an isometry is
either elliptic, or parabolic or hyperbolic, and that if $X$ is an
$\mathbbm{R} -$tree an isometry cannot be parabolic. An elliptic
group cannot contain a hyperbolic or a parabolic element, a
loxodromic group cannot contain a parabolic element. In order to
simplify the notation, if $\phi : G\to Isom(X)$ is an action of the
group $G$, we denote by $gx$ the image of $x$ under the isometry
$\phi(g)$.

\begin{proposition}\label{1.4} Let the pair $(X,G)$ satisfy (C).
Assume that each element of $G$ is either elliptic or parabolic.
Then $G$ is either elliptic or parabolic; if $X$ is an $\mathbbm
R$-tree, and $G$ is finitely generated, then $G$ is elliptic.
\end{proposition}

The proof of this proposition is of dynamical nature and based on the
following

\begin{lemma}\label{1.5}
(See~\cite{coo_del_pap}, chap 9, lemma~\ref{2.3}). Let $X$ be a
$\delta-$hyperbolic space. Le $g$, $h$ be two elliptic or parabolic
isometries of $X$. Suppose that $ \min (d(gx, x), d(hx,x)) \geqslant
2 \langle gx, hx\rangle_x + 6 \delta$. Then $g^{-1}h$ is hyperbolic.
\end{lemma}

Recall that the Gromov product $\langle x,y\rangle_z$ is defined as
$1/2(d(x,z)+d(y,z)-d(x,y))$

\begin{proof}[Proof of Proposition~\ref{1.4}.]
Let us first consider the case where $X$ is an $\mathbbm R$-tree,
which we denote by $T$, and $G$ is finitely generated. Recall that
projection of a point $x$ in a $CAT(0)$ space onto a complete convex
subset $Y$ is the unique closest point to $x$ in $Y$
(see~\cite{bri_hae}, page 176). We claim that in an $\mathbbm
R$-tree $T$, if $g$ is some elliptic isometry, and $T^g$ the subtree
of fixed points of $g$, then for every $x$ the midpoint of the
segment $[ x, g x ]$ is the projection of $x$ on $T^g$: indeed let
$p$ be this projection, so that the image of the segment $[x,p]$ is
$[gx,p]$; if the Gromov product $\langle x,gx\rangle_p=d$ is
strictly positive, we can consider the point $q\in [p,x]$ s.t.
$d(p,q)=d$; it is fixed by $g$ as it belongs to $[p,x]$ and it is
the unique point on this segment with $d(q,p)=d$, but $q$ is closer
than $p$ to $x$, contradiction. Thus $\langle x,gx\rangle_p=0$, and
as the two segments $[x,p]$ and $[gx,p]=g[x,p]$ have the same
length, $p$ is the midpoint of $[x,gx]$. For every subset $\Sigma
\subset G$, let $T_{}^{\Sigma}$ be the fixed subset of $\Sigma$. Let
$\{ g_1, \ldots .g_n \}$ be a finite generating subset of $G$, and
let us prove by induction that $T^{\{ g_1, \ldots g_n \}}$ is not
empty. For $n = 1$ this is the hypothesis. Suppose that $T^{\{ g_1,
\ldots g_{n - 1} \}} \cap T^{\{ g_n \}} = \varnothing$. The minimal
distance between these two subtrees is achieved along a segment $[a,
b]$, with $a \in T^{\{ g_1, \ldots, g_{n - 1} \}}$, $b \in T^{\{ g_n
\}}$. Let $x_0$ be the midpoint of this segment: $x_0 \not\in T^{\{
g_n \}}$. Therefore $b \in [ x_0, g_n x_0 ]$ is the midpoint. As
$x_0 \not\in T^{\{ g_1, \ldots g_{n - 1} \}}$, we have that $x_0
\not\in T^{\{ g_i \}}$, for some $i$. The intersection $T^{\{ g_i
\}} \cap [ a, x_0 ]$ is a segment $[a,c] $; the right extremity $c$
of this segment is the projection of $x_0$ on $T^{\{ g_i \}}$, and
therefore $c \in [ x_0, g_i x_0 ]$ is the midpoint. Thus $x_0 \in [
g_i x_0, g_n x_0 ]$ and, in other words, $\langle g_i x_0, g_n
x_0\rangle_{x_0} = 0$. Lemma~\ref{1.5} applies and proves that the
isometry $g_i g_n$ is hyperbolic, a contradiction.

Suppose now that $X$ is a proper geodesic hyperbolic space. Let $G$
be as in the statement, and $x_0 \in X$ be some base-point. If the
orbit $Gx_0$ is bounded, then it is a bounded $G$ invariant set, and
$G$ is elliptic. Assume that $Gx_0$ is not bounded. We consider the
set $\overline{Gx_0} \cap
\partial X$.

1) Assume that this set has only one point $a$. It must be $G$
invariant. Let us prove that $G$ is parabolic. Suppose that $G$
fixes another point $b$ on the boundary. Then it acts on the union
$Y$ of geodesic lines between $a$ and $b$. Let $L\subset Y$ be a
geodesic between $a$ and $b$, so that every point in $Y$ is at
distance $< 100\delta$ of $L$. Let $x_0\in L$; as $ \overline{Gx_0}
\cap
\partial X =\{a\}$, we can find two isometries $g,h$ in $G$ such that
$d(x_0,g x_0) >1000\delta$, $d(x_0,hx_0)> d(x_0,g x_0)+1000\delta$
and the projections of $gx_0$ and $hx_0$ on $L$ are on the right of
$x_0$.

Considering these projections of $gx_0$ and $hx_0$ on $L$, we see
that $d(x_0,hx_0) \geq d(x_0,gx_0)+d(g x_0,hx_0)- 200\delta$, thus
$\langle x_0,hx_0\rangle_{g(x_0)} \leq 100\delta$. By isometry, we
get $\langle g^{-1}x_0,g^{-1}hx_0\rangle_{x_0} < 100\delta < 1/2
(\min(d(x_0, g^{-1}x_0), d(x_0,g^{-1}hx_0)) -3\delta$, and $h$ must
be hyperbolic by Lemma~\ref{1.5}.

2) Assume that $\overline{Gx_0} \cap
\partial X$ has at least two points, $a, b \in \overline{Gx_0} \cap
\partial X$. There exists two sequences $g_n$ and $h_n$ such that
$g_n x_0 \rightarrow a$, and $h_n x_0 \rightarrow b$. Then $d(g_n
x_0, x_0) \rightarrow \infty$ as well as $d(h_n x_0, x_0)$, but
$\langle g_n x_0, h_n x_0\rangle_{x_0} \rightarrow \langle a,
b\rangle_{x_0}$ and remains bounded (by the very definition of the
Gromov boundary). Lemma~\ref{1.5} applies and we get a
contradiction.
\end{proof}

\begin{corollary}\label{1.6}
Let $(G,X)$ satisfy $(C)$. If $G$ has a subgroup of finite index
which is elliptic or parabolic, the $G$ is also elliptic or
parabolic.
\end{corollary}

\begin{proof}
No element of $G$ can be hyperbolic, as any power of a hyperbolic
element is hyperbolic.
\end{proof}

\begin{proposition}\label{1.7}
Let the pair $(X,G)$ satisfy (C). If $G $ is elliptic, then it has
an orbit of diameter $\leqslant 100 \delta$. If, furthermore, $X$ is
$CAT(0)$, then $G$ has a fixed point.
\end{proposition}

\begin{proof}
In a metric space, the radius of a bounded set $Y$ is the infimum of
$r$ s.t. there exists a $x$ with $Y \subset B(x, r)$. A center is a
point $c$ s.t. $Y \subset B(c, r')$ for every $r' > radius(Y)$. The
proof of the proposition~\ref{1.7} is a direct consequence of the
following generalization of Elie Cartan center's theorem
\cite{bri_hae}, II.2.7.
\end{proof}

\begin{proposition}\label{1.8}
In a proper geodesic $\delta$-hyperbolic space, the diameter of the
set of the centers of a bounded set is $\leqslant 100 \delta$. In a
complete $CAT(0)$ space, every bounded set admits a unique center.
\end{proposition}

\begin{proof}
The second point is proved in~\cite{bri_hae}, II.2.7. Let us prove
the first assertion. Let $a, b$ be two centers and suppose that
$d(a, b) > 100 \delta$. Let $c$ be a midpoint of $a, b$. Let us
prove that for every $x$ in $Y$, $d(y, c) < r - 10 \delta$ and in
such way get a contradiction. By assumption $d(a, x)$ and $d(b, x)$
are less than $r + \delta$. By the 4 points definition of
$\delta$-hyperbolicity (\cite{coo_del_pap} prop.~1.6) we know that $
d(x,c)+d(a,b) \leq \max( d(x,a)+d(b,c), d(x,b)+d(x,c)) - 2\delta$.
 As $d(b,c)=d(a,c)=1/2d(a,b)>50 \delta$ we get that $ d(x,c) \leq
 \max(d(x,a), d(x,b)) -48 \delta \le r-48 \delta$ and we are done.
\end{proof}

\begin{proposition}\label{1.9} Let the pair $(X,G)$ satisfy (C).
If the $G$-orbit of some point of $\partial X$ is finite and has at
least $3$ elements, then $G$ is elliptic.
\end{proposition}

\begin{proof}
If the orbit is finite and has at least $3$ elements $w_1, \ldots
w_k$, let us construct a bounded orbit of $G$ in $X$. For every
triple of different points $w_i, w_j, w_k$ in this orbit, let us
consider the set $C_{ijk}$ consisting of all points being at a
distance less than $24 \delta$ from all geodesics between $w_i$,
$w_j$, and $w_k$. By hyperbolicity this set is not empty and has
diameter $\leq 100 \delta$. This follows from~\cite{coo_del_pap},
chap 2 prop.2.2, page 20. A finite union of bounded sets is bounded.
Therefore, the union of the sets $C_{ijk}$ is a bounded $G$
invariant set.
\end{proof}

\begin{proposition}\label{1.10}
Let $X$ be a $\delta -$hyperbolic graph of bounded valence. If $G
\subset \Isom{X}$ is loxodromic, then there exists an epimorphism
$m: G \rightarrow \mathbbm{Z}$ or $m: G \rightarrow
\mathbbm{D}_{\infty}$ such that $\ker m$ is elliptic.
\end{proposition}

\begin{proof}
We will construct a combinatorial analogue of the Busemann cocycle
(compare~\cite{ripssela}). As $G$ is loxodromic, the action of $G$
fixes two points $w^{\pm}$ at infinity. It contains a subgroup of
index at most two $G^+$ which preserves these two points, and
contains some hyperbolic element $h$. Let $U$ be the union of all
geodesics between these two points at infinity, and choose a
preferred oriented line $L$ between this two points. If $x\in U$,
there exist a point in $L$ such that $d(x,p(x)) < 24\delta$
(\cite{coo_del_pap}, chap 2 prop.2.2, page 20). Choose such a point
and call it a projection of $x$. If $x \in U $, let $R(x) = \{y \in
U | d(x, y) > 1000 \delta$, and the projection of $y$ to $L$ is on
the right to that of $x$ $\}$. Note that our hypothesis implies that
for every pair $x, y$, $ \{ R(y) / R(x) \}$ is contained in the ball
centered at $y$ and of radius $d(x,y)+2000\delta$ and is therefore
finite: by definition, a point of $R(y)$ which is at distance $>
d(x,y)+2000\delta $ from $y$ must project on $L$ on a point which is
at the distance $> 1000 \delta $ of $x$. Note also that if $h$ is
hyperbolic, $R(h^n x)$ is strictly contained in $R(x)$ if $n$ is
$\gg 1.$ Let $c(x, y) =Card \{ R(y) \setminus R(x) \} - Card \{ R(x)
\setminus R(x) \}$. Note that $c(y, x) + c(x, y) = 0$, and that
$c(x, y) + c(y, z) = c(x, z)$. Moreover, if $g$ is in $G^+$ then
$R(gx) = gR(x)$. Choose some point $x_0\in U$. The formula $m(g) =
c(x_0, gx_0)$ defines a non trivial morphism $G^+ \rightarrow
\mathbbm{Z}$. The orbit of $x_0$ under the action of the kernel of
$m$ is bounded, contained in $B(x_0, 2000 \delta)$, and $\ker m$ is
elliptic. If $G / G^+$ is not trivial, and $\varepsilon \in G
\setminus G^+$ then $m - m (\varepsilon g \varepsilon^{- 1}) = d ( g
)$ extends to a non trivial epimorphism
 $G \rightarrow \mathbbm{D}_{\infty}$.
\end{proof}

\begin{proof}[Proof of Theorem~\ref{1.1}.]

Let $H_1$ be the rigid stabilizer of the first level of $G$. It is a
product of $n$ subgroups of $G$, $H_1 = L_1 \times \ldots \times
L_n$ conjugate in $G$.

i) Suppose first that $L_1$ contains no hyperbolic element.

Then $L_1$ has either (1) a bounded orbit or (2) a unique fixed
point $w$ at infinity.

(1) In the first case, let $C_1 = \{ x| \forall g \in L_1, d(gx, x)
< 100 \delta \}$ (by Proposition 2 this set is nonempty). As $L_2$
commutes with $L_1$ it preserves $C_1$. Being conjugate to $L_1$,
every orbit of $L_2$ is bounded. If $x_0 \in C_1$ and $D=diam(L_2
x_0)$, we see that the diameter of $(L_1 \times L_2) x_0$ is
$\leqslant D + 2\cdot 100 \delta$, hence $L_1 \times L_2$ is
elliptic, and the set $C_2 = \{ x| \forall g \in L_1 \times L_2,
d(gx, x) \leq 100 \delta \}$ is not empty (Proposition 2). By
induction we prove that $C_k = \{ x| \forall g \in L_1 \times L_2
\times \ldots \times L_k, d(gx, x) < 100 \delta \}$ is not empty;
thus $G$ admits a subgroup of finite index which is elliptic, and
$G$ is itself elliptic.

(2) In the second case, the unique fixed point $w$ is stable under
the action of the subgroup $L_2 \times \ldots \times L_n$, and $G$
has a subgroup of finite index which is parabolic, thus $G$ is
parabolic itself.

ii) Suppose $L_1$ contains some hyperbolic element $h$. Let
$w^{\pm}$ be the two distinct fixed points of $h$ at infinity. As
$L_2 \times \ldots \times L_n$ commutes with $h$ this group fixes
this set. Now $L_2$ contains a hyperbolic element $h_2$, conjugated
to $h$: thus $h_2$ has the same fixed points at infinity as $h$, and
$H_1$ must also fix the set $\{w,w^-\}$. Thus the orbit of $w^{\pm}$
is finite and Proposition ~\ref{1.9} applies. The orbit of $G$
cannot have more than 2 elements unless $G$ is elliptic: therefore
it has exactly two elements, and $G$ is loxodromic.This proves a).

To prove b) apply Proposition~\ref{1.10}. Proposition~\ref{1.7} and
Proposition~\ref{1.8} (the unique center for bounded sets) give rise
to the desired fixed point for c). Claim d) follows from the fact
that between two points at infinity in a CAT(-1) space there exists
a unique geodesic (visibility property). For claim e), let $w$ be an
end of a tree fixed by the group $G$. Let $t\to r(t)$ be a geodesic
ray converging to $w$. Note that when the point $x$ is fixed, the
function $t\to d(x,r(t))-t$ is constant for $t>>1$. The value of the
constant $b_w(x)$ is called the Busemann function associated to $w$
(see~\cite{bri_hae}, Chap II.8. for a study of Busemann functions in
CAT(0) spaces). If the point $w$ is fixed by some isometry $g$, then
$g.r(t)$ is another ray converging to $w$. But two rays converging
to the same point in a tree must coincide outside a compact set.
Therefore $d(gr(t),r(t))=b(g)$ is constant for $t>>1$, and this
constant is $b_{g.w}-b_w$. By construction, $g\to b(g)$ is a
homomorphism from $G$ to $\mathbbm{R}$, which is non-trivial unless
every element of $G$ is elliptic, and takes values in $\mathbbm{Z}$
if $X$ is a combinatorial tree. Suppose that the restriction of $b$
to $L_1$ is trivial. Then $L_1$ consists of elliptic elements. Since
$G$ is finitely generated, $L_1$ is finitely generated as well. Thus
$L_1$ is elliptic and i) applies. Otherwise, $L_1$ contains a
hyperbolic element and ii) applies.


Theorem~\ref{1.1} is proved.
\end{proof}

\section{An indicable branch group}
\label{sec:example1}

Let $\mathcal{G}$ be a branch group acting on a rooted tree $\Tree$.
It is proved in~\cite{gri:branch} that, if $N \triangleleft
\mathcal{G}$ is a nontrivial normal subgroup, then the group $N$
contains the commutator subgroup of the rigid stabilizer
$rist_{\mathcal{G}}(n)'$, for some level $n$. As
$rist_{\mathcal{G}}(n)$ is of finite index in $G$,
$G/rist_{\mathcal{G}}(n)$ is finite, $G/rist_{\mathcal{G}}(n)'$ is
virtually abelian and we have:

\begin{proposition}\label{2.1}
 A proper quotient $\mathcal{G}/N$ of a branch group is a virtually
abelian group.
\end{proposition}

We construct in this section an example of a finitely generated
branch group which surjects onto the infinite cyclic group. The
construction starts from the finitely generated torsion 2-group
firstly defined in~\cite{grigorch:80} and later studied
in~\cite{Grigorch:84} and other papers (see also the Chapter VIII of
the book~\cite{harpe}).

We will list briefly some properties of $G$ that will be used later.

Let $(\mathcal{T}, \emptyset)$ be the rooted binary tree whose
vertices are the finite sequences of $0, 1$ with its natural tree
structure (see~\cite{harpe}, VIII.A for details), the empty sequence
$\emptyset$ being the root. If $v$ is a vertex of $\mathcal{T}$ we
denote by $\mathcal{T}_v$ the subtree consisting of the sequences
starting in $v$. In other words, the subtree $\mathcal{T}_v$ of
$\mathcal{T}$ consists of vertices $w$ that contain $v$ as a prefix.
Deleting the first $|v|$ letters of the sequences in $\mathcal{T}_v$
yields a bijection between $\mathcal{T}_v $ and $\mathcal{T}$,
called the canonical identification of these trees.

\begin{figure}[h]
\begin{center}
\includegraphics{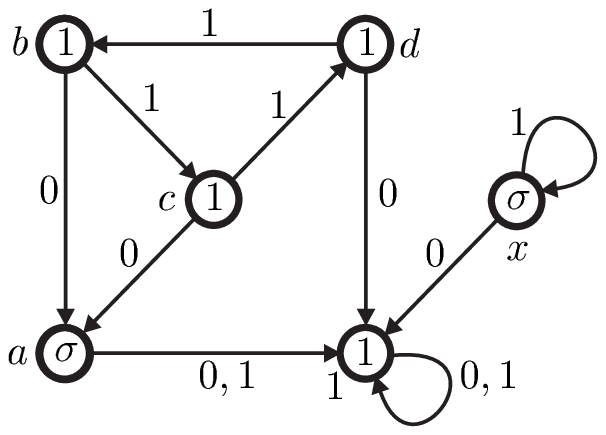}
\caption{The automaton defining $L$} \label{D1}
\end{center}
\end{figure}

The group $G$ (see~\cite{harpe}, VIII.B.9 for details) acts
faithfully on the binary rooted tree $(\Tree, \emptyset)$ and is
generated by four automorphisms $a, b, c, d$ of the tree where $a$
is the rooted automorphism permuting the vertices of the first
level, while $b, c, d$ are given by the recursive rules
\[
 b = (a, c), c = (a, d), d = (1, b).
\]
This means that $b$ does not act on the first level of the tree, it
acts on the left subtree $\mathcal{T}_0$ as $a$ and acts on the
right subtree $\mathcal{T}_1$ as $c$, with similarly meaning of the
relations for $c$ and $d$. Here we use the canonical identifications
of $\Tree$ with $\Tree _i, i=0,1$. An alternative description of $G$
is that it is the group generated by the states of the automaton
drawn on the figure~\ref{D1}.

The group $G$ is $3$-generated as we have the relations
\[a^2 = b^2 = c^2 = d^2 = bcd = 1\]
 there are many other relations and $G$ is not finitely presented
\cite{Grigorch:84}.

In order to study groups acting on the binary rooted tree $\mathcal
T$, it is convenient to use the embedding
\[
 \psi: \text{$\mathop{Aut} (\mathcal{T})$} \longmapsto
 \text{$ \mathop{Aut} (\mathcal{T})$} \wr \Sym{2},
\]
\[ g \longmapsto (g_0, g_1) \alpha .\]
In this description $\Sym{2}$ is a symmetric group of order 2,
 $\alpha \in \Sym{2}$ describes the action of $g$ on the first level of
 the tree and the sections $g_0, g_1$ describe the action of $g$ on the of subtrees $\mathcal{T}_0$ $\mathcal{T}_1$.
 We will usually identify the element $g$ and its image $(g_0, g_1) \alpha$. Relations of
 this type will be often used below.

Let $x$ be the automorphism of $\mathcal{T}$ defined by the
recursive relation $x = (1, x)a$. This automorphism is called the
adding machine as it imitates the adding of a unit in the ring of
diadic integers~\cite{gns}. An important property of $x$ is that it
acts transitively on each level of $\Tree$ and therefore has
infinite order.

Let $L = \langle x, G\rangle$ be the subgroup of $\mathop{Aut}
(\mathcal{T})$ generated by $G$ and the adding machine $x$.

\begin{theorem}\label{2.2}
The group $L$ is branch, amenable, and has infinite abelization.
\end{theorem}

The next two lemmas are the first steps towards the proof of the
fact that $L$ is a branch group.

\begin{lemma}\label{2.3}
The following formulas hold in the group L:

$[x, a] = (x^{- 1}, x)$,

$[x, d] = (x^{- 1} b x, b)$

$[[x, a], d] = (1, [x, b])$,

$(1, [[x, b ], c]) = [[[x, a], d], b]$.
\end{lemma}

\begin{proof} This follows by direct computation:

 $[x, a] = x^{- 1} axa = a (1, x^{- 1}) a (1, x) aa = (x^{- 1}, 1) (1, x) = (x^{- 1}, x)$.

 $[x, d] = x^{- 1} dxd = a (1, x^{- 1}) (1, b) (1, x) a (1, b) =
 (x^{- 1} b x, b)$

 $[[ x, a], d] = [(x^{- 1}, x), (1, b)] = ([x^{-1}, 1], [x, b]) = (1, [x, b])$.

 $[[[x, a], d], b] = [(1, [x, b]), (a,
 c)] = (1, [[x, b], c])$.
\end{proof}

\begin{lemma}\label{2.5}
The group L is self-replicating, and hence level transitive.
\end{lemma}

\begin{proof}
Consider the elements $b = (a, c)$, $c = (a, d)$, $d = (1, b)$, $aba
= (c, a)$, $xa = (1, x)$. They stabilize the two vertices of the
first level of $\mathcal{T}$, and their projections on $\mathop{Aut}
(\mathcal{T}_1) \simeq \mathop{Aut} (\mathcal{T})$ are $c$, $d$,
$b$, $a$, $x$, i.e. the generators of $L$. Note that these elements
generate $L$. Hence the projection of $s t_L(1)$ on $\mathop{Aut}
(\mathcal{T}_1)$ is $L$ modulo the canonical identification of
$\mathcal{T}$ and $\mathcal{T}_1$. The conjugation by $a$ permutes
the coordinates of elements in $st_L(1)$, hence the same holds for
the first projection. The self-replicating property (definition 4)
follows by induction on the level. The level transitivity is an
immediate consequence of the transitivity of $L$ on the first level
and the self-replicating property.
\end{proof}

Let
\[K = \langle [a, b] \rangle^G, S = [ \langle x \rangle, G ]^L,\]
\[R = \langle K, S, \gamma_2(L) \rangle^L=KS\gamma_2(L).\]
These subgroups will play an important role in our further
considerations.

\begin{lemma}\label{2.7}
We have the following inclusions: $\gamma_2(G) \succeq \gamma_2(G)
\times \gamma_2(G), K \succeq K \times K,$and $R \succeq S \times
S$.
\end{lemma}

\begin{proof}
The first two inclusions are known~\cite{gri:hilb, gri:branch}.

Using the commutator relations and the fact that conjugation by $a$
permutes the coordinates we have
\[(1, [c, x]) = [(a, c), (1, x)] = [b, xa] = [b, a] [b, x] [[b, x], a] \in R, \]
\[(1, [x, b]) = [(x^{- 1}, x), (1, b)] = [[x, a], d] \in
R,\] by Lemma~\ref{2.3},
\[(1, [a, x]) = a [(a, c), (x, 1)] a = a [b, (x, 1)] a =
ab^{- 1} (x^{- 1}, 1) b (x, 1) a.\] But $x = (1, x) a$ and $axa =
(x, 1)a$ which leads to
 \[(1, [a, x]) = ab^{- 1} x^{- 1} abaxa = ab^{- 1} ab [b, axa] = [a, b] [b, axa].\]

Now we have
 \[[ b, axa ] = a [ aba, x ] a \in S,\] and $[ a, b ]
\in K$ which gives
 $(1, [ a, x]) \in R.$

Finally
\[(1, [x, d]) = (1, [x, bc]) = (1, [x, c] [x, b] [[x,
b], c]) = (1, [x, c]) (1, [x, b]) (1, [[x, b], c])\] and
\[( 1, [ [ x, b ], c ] ) = [ [ [ x, a ], d ], b ] \in R\] by
Lemma~\ref{2.3}.
Therefore the elements $( 1, [ a, x ] )$, $( 1, [b, x ] )$, $( 1, [ c, x ] )$,
$( 1, [ d, x ] )$ belong to $R$ and, as
$S = \langle [ a, x ], [ b, x ], [ c, x ], [ d, x ] \rangle^L$, the lemma is proved.
\end{proof}

\begin{lemma}\label{2.8}
We have the inclusion: $\gamma_2(L) \succeq \gamma_2(L) \times
\gamma_2(L)$.
\end{lemma}

\begin{proof}
Consider the subgroup $Q = \langle d, c, aca, xa \rangle \subset L$.
As $d = (1, b)$, $c = (a, d)$, $aca = (d, a)$, $xa = (1, x)$, the
group $Q$ is a subdirect product in $D_4 \times L$ where $D_4\simeq
\langle a, d \rangle $ is a dihedral group of order 8. As
$\gamma_2(D_4) = 1$ we get
\[\gamma_2(Q) = (1, \gamma_2(L)),\]
\[\gamma_2(aQa)=(\gamma_2(L),1),\]
and therefore
\[\gamma_2(L) \succeq \gamma_2(L) \times \gamma_2(L).\]
\end{proof}

\begin{lemma}\label{2.9}
The group $L$ is a weakly regular branch group over $R$.
\end{lemma}

\begin{proof}
We know that $K \succeq K \times K, \ \gamma_2(G) \succeq
\gamma_2(G) \times \gamma_2(G)$, $\gamma_2(L) \succeq \gamma_2(L)
\times \gamma_2(L)$ and $R \succeq S \times S$. But $R$ is generated
by $S,$ $\gamma_2(L)$ and $K$. This implies the statement.
\end{proof}

In order to prove that $L$ is a branch group, we consider its
subgroup $P = \langle R, \langle x^4 \rangle \rangle^L$.

\begin{lemma}\label{l:2.10}
The group $P$ has finite index in $L$.
\end{lemma}

\begin{proof}
Every element $g \in L$ can be written as a product $g = x^i a^j c^k
d^l hfx^{4 t}$, where $h \in [ G, G ], f \in S, i \in \{ 0, 1, 2, 3
\}$, $j, k, l \in 0, 1, t \in \mathbb{Z}$. This implies that the
index of $P$ in $L$ is $\leqslant 128$.
\end{proof}

Let $P_n \simeq P \times \cdots \times P \subset \mathop{Aut}
(\mathcal{T})$ ($2^n$ factors) be the subgroup of $\mathop{Aut}
(\mathcal{T})$ that is the product of $2^n$ groups isomorphic to $P$
that act on the corresponding $2^n$ subtrees rooted at the vertices
on the $n$-th level.

\begin{lemma}\label{2.11}
The group $L$ contains $P_n$ for every $n$.
\end{lemma}

\begin{proof}
For $n = 0$ the statement is obvious. For $n = 1$, let us consider
$(xa)^4 = (1, x^4)$ which is an element of $L$. As $L$ is
self-replicating, for any given element $h \in L$ there exists an
element $k$ in $L$ s.t. $k = (f, h)$. Conjugating $(1, x^4)$ by an
element of $L$ of the form $(f, h)$, we get that $(1, (x^4)^h) \in
L$. But $P$ is generated by conjugates of $x^4$. This together with
Lemma~\ref{2.9} proves the inclusion $1 \times P \succeq L$. The
inclusion $P \times 1 \succeq L$ is obtained by conjugating $L$ by
$a$. Then we get that $P \times P = P_1 < L$.

In order to prove the lemma for $n = 2$ we observe that
\[L \ni [x,a] = (xa)^2 = (1, x^2) = (1, 1, x, x)_2\]
(the index $2$ indicates that we rewrite the considered while
considering its action on the second level; we will use such type of
notations for further levels as well). Multiplying $(1, x^2)$ (which
is in $L$) by
\[(1, [ x, a ]) = (1, 1, x^{- 1}, x)_2,\]
we get $(1, 1, 1, x^2)_2 \in L$. Therefore $(1, 1, 1, x^4)_2 \in L$
and hence $P_2 < L$ (by level transitivity and the self-replicating
property of $L$ we see that $(x^4,1,1,1), (1,x^4,1,1)$ and
$(1,1,x^4,1)$ also belong to $L$.

Let us prove the lemma by induction on $n \geq 2$. Suppose that, for
every $k \leqslant n$, the inclusion $P_k < L$ holds and let us
prove that $P_ {n + 1} < L$. Consider the element $\mu$
\[L \ni \mu = (1, \ldots, 1, x^4)_{n - 2} = (1, \ldots, 1, x^2, x^2)_{n - 1,}\]
As $L$ is self-replicating, there exists an element $\rho \in
St_G(u_{n - 2})$, where $u_{n - 2}$ is the last vertex on the $(n -
2)$-th level, whose projection at this vertex is equal to $b$. We
have

\begin{equation}
\begin{array}{l}
 L \ni [ \mu, \rho ] = (1, \ldots, 1, [ x^2, a ], [ x^2, c ])_{n - 1}\\
 = ( 1, \ldots, 1, [ x, a ], [ x, d ] )_n
 = ( 1, \ldots, 1, x^{- 1}, x, x^{- 1} bx, b )_{n + 1},
\end{array}
\end{equation}
As $b^2=1$ we get the relation
\[ [ \mu, \rho ]^2 = ( 1, \ldots, 1, x^{- 2}, x^2, 1, 1 )_{n + 1} \]
 Now we have
\begin{equation}
\begin{array}{l}
 L \ni \eta = ( 1,,, 1, x^4, 1 )_n = ( 1, \ldots, 1, x^2, x^2, 1, 1 )_{n
 + 1},\\
 \left[ \mu, \rho \right]^2 \eta = ( 1, \ldots, 1, x^4, 1, 1 )_{n + 1}
\end{array}
\end{equation}
and we come to the conclusion that $1 \times 1 \times \ldots \times
1 \times P \times 1 \times 1 \succeq L$, hence
$\underbrace{P\times \ldots \times P}_{2^{n+1}} \succeq
R$, and $P_{n + 1} < L$, as $L$ is level transitive.
\end{proof}

We can now prove that $L$ is a branch group. This group acts
transitively on each level of the rooted tree $\mathcal {T}$, and
contains $P_n$ for every $n = 1, 2 \ldots$~. In order to prove that
it is branch, as $P_n < rist_L(n)$, and $L$ is level transitive, it
is enough to check that $P_n$ has finite index in $L$. We have the
following diagram
\[
\begin{array}{lllllll}
 L & & & & & & \\
 \uparrow & \psi_n & & & & & \\
 st_L(n) & \twoheadrightarrow & \tilde{H} & < & L & \times \ldots
 \times & L\\
 \uparrow & & & & & & \\
 rist_L(n) & & \uparrow & & \uparrow & \ldots & \uparrow\\
 \uparrow & \psi_n & & & & & \\
 P_n & \twoheadrightarrow & \tilde{P}_n & = & P & \times \ldots \times & P
\end{array}
\]
(the vertical arrows are inclusions, $\tilde{H}$ and $\tilde{P}_n$
are $\psi_n$ images of $st_L(n)$ and $P_n$ respectively, where
$\psi_n$ is the $n$-th iteration of $\psi$).

As the group $P$ has finite index in $L$, we get that $\tilde{P}_n$
has finite index in $\tilde{H}$ and therefore $P_n$ has finite index
in $st_L(n)$ and hence in $L$. This establishes the first statement
of Theorem~\ref{2.2}.

The group $L$ is the self-similar group generated by the states of
the automaton in Figure~\ref{D1}. The diagram of this automaton
satisfies the condition of Proposition 3.9.9 of \cite{nekrash:book}:
it is therefore a bounded automaton in the sense of
Sidki~\cite{sidki:bounded}. This proposition states that an
automaton is bounded if and only if its Moore diagram has the
following property: every two nontrivial cycles are disjoint and are
not connected by a directed path ; a cycle is called trivial if all
of its states represents the identity automorphism of the tree.

It is easy to see that automaton determining the group $L$ satisfies
this property.

By a theorem of Bondarenko and Nekrashevych (Theorem 3.9.12
in~\cite{nekrash:book}) every group generated by the states of a
bounded automaton is contracting. Moreover, by a theorem of
Bartholdi, Kaimanovich, Nekrashevych and Virag~\cite{bknv:amenab}
such a group is amenable. This establishes the amenability of $L$,
as well as its contracting property.

In order to compute the abelianization of $L$, we need to combine
the contracting property of $L$ with a rewriting process which
corresponds to the embedding $\psi$. The combination of this
rewriting process and the contraction property will produce an
algorithm for solving the word problem in $L$: the branch algorithm.
This type of algorithm appeared in~\cite{Grigorch:84} for the first
time: it is a general fact that the branch algorithm solves the word
problem for contracting groups~\cite{savchuk}.

The group
\[
 \Gamma=\langle a,b,c,d,x:a^2=b^2=c^2=d^2= bcd=1 \rangle,
\]
defined by generators and relations, naturally covers $L$. It is
isomorphic to the free product
\[
\mathbb{Z}/2\mathbb{Z} \ast (\mathbb{Z}/2\mathbb{Z}\times
\mathbb{Z}/2 \mathbb{Z}) \ast \mathbb{Z}.
\]
Therefore, the elements in $\Gamma$ are uniquely represented by
words $w=w(a,b,c,d,x)$ in the reduced form (for this free product
structure).

Similarly the group $G$ is naturally covered by the group
\[\langle a,b,c,d: a^2=b^2=c^2=d^2= bcd=1\rangle \simeq \mathbb{Z}/2\mathbb{Z}
\ast ( \mathbb{Z}/2\mathbb{Z}\times \mathbb{Z}/2\mathbb{Z}).\] The
elements in $G$ can be represented by reduced words (with respect to
this free product structure).

Let $w$ be a word representing an element of $\Gamma$, $w = u_1
x^{i_1} u_2 x^{i_2} \ldots u_k x^{i_k} u_{k + 1}$, where $u_i$ are
reduced words in $a, b, c, d$, $u_i$ is nonempty for $i \neq 1, k +
1$, and $i_j \neq 0$, for $j=1,...,k$.

Let us consider the following rewriting process:

1) in each word $u_i$ replace $b, c, d$ by the corresponding element
of the wreath product $L \wr \Sym{2},$ using the defining relations
$b = (a, c), c = (a, d)$ $d = (1, b)$,
 $x=(1, x)a$.

 2) Move all the letters $a$ to right using the relations $a (v_0, v_1)
= (v_1, v_0) a$. Use the relation $a^2 = 1$ for simplification of
words, and take the componentwise product of all involved pairs. One
obtains in such a way a relation of the form
$w=(w_0,w_1)a^{\varepsilon}$ with $\varepsilon \in \{0, 1\}$, which
holds in $L$.

3) Reduce the words $w_i$ in $\Gamma$, obtaining a pair $(w_0,w_1)$
of reduced words.

Note that the length of $w_i, i=0, 1$ is strictly shorter than of
$w$ if at least one letter $a$ appears in the word $w$.

We can represent this rewriting process as a pair $\varphi =
(\varphi_0, \varphi_1)$ (or a product $ \varphi_0 \times \varphi_1$)
of two rewritings $w\to w_0$ and $w \to w_1$. We will apply these
maps to words with an even number of occurences of $a$, i.e. words
representing the elements in $\mathop{st_L}(1)$: in this case
$\varepsilon = 0.$ We can therefore iterate this rewriting procedure
$\varphi$ $n$ times for words representing elements in
$\mathop{st_L}(n)$, and get $2^n$ words $w_{ i_1, \ldots i_n }$ with
$i_j \in \{ 0, 1 \}$. (For formal definition of $\varphi_0,
\varphi_1 $ in case of the group $G$ see~\cite{gri:example}, for $L$
the formal description is similar).

\begin{proposition}\label{2.12}
The rewriting process is 3-step contracting with core $\mathcal{N}$
$=$ $\{1$, $b$, $c$, $d$, $x$, $x^{-1}$, $bx$, $cx$, $dx$, $x^{-1}b$,
$x^{-1}c$, $x^{-1}d$, $x^{- 1} bx$, $x^{- 1} cx$, $x^{- 1}dx\}$.
 In other words: for every word $w$ representing an element in
$\mathop{stab}_L(3)$, $\varphi^3(w)$ consists of 8 words $w_{i, j,
k}, i,j,k \in\{0,1\} $ of strictly shorter length than $w$.
\end{proposition}

\begin{proof}
Let the word $w = w_1 x^{i_1} w_2 x^{i_2} \ldots w_k x^{i_k} w_{k +
1}$ be
 as above and represents an element in $stab_L(3)$. As we already have
noted, if the letter $a$ occurs in some of the $w_i$ then rewriting
process is strictly shortening in one step. In order to study
reduced words without the letter $a$, we will make use of the
relations in Table~\ref{t:relations}.
\begin{table}
\begin{eqnarray*}
& bx=(a,c)(1,x)a=(a,cx)a\\
& cx=(a,d)(1,x)a=(a,dx)a\\
& dx=(1,b)(1,x)a=(1,bx)a\\
& x^{-1}b=a(1,x^{-1})(a,c)=(x^{-1}c,a)a\\
& x^{-1}c=a(1,x^{-1})(a,d)=(x^{-1}d,a)a\\
& x^{-1}d=a(1,x^{-1})(1,b)=(x^{-1}b,1)a\\
& xb=(1,x)a (a,c)=(c,xa)a & B\\
& xc=(1,x)a(a,d)=(d,xa)a &B \\
& xd=(1,x)a(1,b)=(b,x)a& A\\
&bx^{-1}=(a,c)a(1,x^{-1})=(ax^{-1},c) & B\\
&cx^{-1}=(a,d)a(1,x^{-1})=(ax^{-1},d)&B\\
&dx^{-1}=(1,b)a(1,x^{-1})=(x^{-1},b)&A\\
& xbx = (1, x) a (a, c) (1, x) a = (1, x) (c, a) ( x, 1 ) = (
cx, xa ), & A\\
& x^{- 1} bx^{- 1} = a ( 1, x^{- 1} ) (a, c) a ( 1, x^{- 1} ) = ( x^{-
1} c, ax^{- 1} ) & A\\
& xbx^{- 1} = (1, x) a (a, c) a ( 1, x^{- 1} ) = (1, x) (c, a) (
1, x^{- 1} ) = ( c, xax^{- 1} ) & B\\
& x^{- 1} bx = a ( 1, x^{- 1} ) (a, c) (1, x) a = ( x^{- 1} cx, a ),
&
C\\
& xcx = (1, x) a (a, d) (1, x) a = (1, x) ( d, a ) ( x, 1 ) = (
dx, xa ) & A\\
& x^{- 1} cx^{- 1} = a ( 1, x^{- 1} ) (a, d) a ( 1, x^{- 1} ) = ( x^{-
1}, 1 ) ( d, a ) ( 1, x^{- 1} ) = ( x^{- 1} d, ax^{- 1} ) & A\\
& xcx^{- 1} = (1, x) a (a, d) a ( 1, x^{- 1} ) = (1, x) ( d, a ) (
x^{- 1}, 1 ) = ( dx^{- 1}, xa ) & A\\
& x^{- 1} cx = a ( 1, x^{- 1} ) (a, d) (1, x) a = ( x^{- 1} dx, a )
&
C\\
& xdx = (1, x) a (1, b) (1, x) a = (1, x) ( b, 1 ) ( x, 1 ) = (
bx, x ) & A\\
& x^{- 1} dx^{- 1} = a ( 1, x^{- 1} ) (1, b) a ( 1, x^{- 1} ) = ( x^{-
1}, 1 ) ( b, 1 ) ( 1, x^{- 1} ) = ( x^{- 1} b, x^{- 1} ) & A\\
& xdx^{- 1} = (1, x) a (1, b) a ( 1, x^{- 1} ) = ( b, 1 ) & A\\
& x^{- 1} dx = a ( 1, x^{- 1} ) (1, b) (1, x) a = ( x^{- 1} bx, 1 )
&
C\\
& x^2 = ( x, x ) & A\\
& x^{- 2} = ( x^{- 1}, x^{- 1} ). & A
\end{eqnarray*}
\caption{Some relations in $L$}\label{t:relations}
\end{table}

Observe that $w$ is a product of subwords in the form presented by
the left side in the relations in Table~\ref{t:relations}, followed
by an element of the set
\[
\{1,b,c,d,xbx,cx,dx,xb,xc,xd,x^{-1}b,x^{-1}c,x^{-1}d,bx^{-1},cx^{-1},dx^{-1}\}.
\] In all relations marked by $A$ or $B$ the rewriting process gives
shortening in one step (case $A$) or in two steps (case $B$); in the
latter case note the presence of the letter $a$, which insures
reduction of length in one more step.

If the word $w$ is not shortened after applying twice the rewriting
procedure, then either it belongs to $N$, or it is of the form $\ast
x^{- 1} \ast x^{} \ldots x^{- 1} \ast x \ast t$, with $\ast \in \{
b, c, d \}$except for the the first or last $\ast$ which may
 also represent the unit, and $t\in \{x^{-1}b,x^{-1}c,x^{-1}d,bx^{-1},cx^{-1},dx^{-1}\}$.

Let $x^{- 1} bx = \tilde{b}$, $x^{- 1} cx = \tilde{c}$, $x^{- 1} dx
= \tilde{d}$. These elements are of order two, and satisfy the
relations $\tilde{b} = (\tilde{c}, a), \tilde{c} = (\tilde{d}, a),
\tilde{d} = (\tilde{b}, 1)$. Since these relations are of the same
form as the relations that hold for $b$, $C$ and $d$, the group
$\tilde{G}$ generated by $\langle a, \tilde{b}, \tilde{c}, \tilde{d}
\rangle$ is isomorphic to $G$.

Let $A < L$ be the subgroup generated by $\langle b, c, d,
\tilde{b}, \tilde{c}, \tilde{d} \rangle$. Note that $A$ stabilizes
the first level of the tree. Consider the embedding $\psi: A
\rightarrow \tilde{G} \times G$ obtained by projecting the elements
of $A$ on the left and right subtrees (we use the same notation
$\psi$ for the embedding as before).

\begin{lemma}\label{2.13}
The group $\psi(A)$ is a subdirect product of finite index in
$\tilde{G} \times G$.
\end{lemma}

\begin{proof}
We have $\tilde{c} = (\tilde{d}, a)$, $c = (a, d)$, $\tilde{d} =
(\tilde{b}, 1)$, $d = (1, b)$, $\tilde{b} \tilde{d} = \tilde{c}$ and
$bd = c$. Therefore the projection of $\psi(A)$ on each of two
factors is onto.

Let $B = \langle b \rangle^G$ and $\tilde{B} = \langle \tilde{b} \rangle^{\tilde{G}}$.

As $d = (1, b) \in \psi(A)$ and as for every $g \in G$ there exists
some $h$ s.t. $(g, h) \in \psi(A)$, we see that $(1, gbg^{-1}) \in
\psi(M)$. Therefore the group $1 \times B$ is contained in $\psi(A)$
and, by a symmetric argument, $\tilde{B} \times 1$ is contained in
$\psi(A)$. Thus $\tilde{B} \times B < \psi(A) < \tilde{G} \times G$.
But the groups $B$ and $\tilde{B}$ have finite index in $G$ and
$\tilde{G}$, respectively, and the lemma is proved.
\end{proof}

We now finish the proof of Proposition~\ref{2.12}. Consider a
reduced word $u$ which represents an element of $L$. Suppose that
this element stabilizes the first level but is not shortened after
applying twice the rewriting process. The word $u$ has to be of the
form $u=wb$, where $w$ represents an element of $A$ and $t\in
\{x^{-1}b,x^{-1}c,x^{-1}d, bx^{-1},cc^{-1}, d{x^{-1}}\}$. Rewrite it
as a word in the letters $\langle b, c, d, \tilde{b}, \tilde{c},
\tilde{d} \rangle$. Use the relations $\tilde{b} = (\tilde{c}, a),
\tilde{c} = ( \tilde{d}, a ), \tilde{d} = ( \tilde{b}, 1 )$, $b =
(a, c), c = (a, d)$ $d = ( 1, b )$ to rewrite it as an element $(
\tilde{w}_0, w_1 )$ of $\tilde{G} \times G$.

 Recall that, endowed with its natural system of
generators, the group $G$ is one step contracting with core
$\mathcal{N}_0 = \{ 1, b, c, d \}$ (and contracting coefficient
$\frac{1}{2}$~\cite{Grigorch:84}). In other words,applying the
rewriting procedure to reduce a word $v$ in $a, b, c, d$ with an
even number of occurrences of the letter $a$ yields a couple a words
of length $\leqslant 1 / 2|v|$ unless $v \in \{ 1, b, c, d \}$. More
precisely, if $v \to (v_0,v_1)$ is obtained by rewriting in the
group $G$, then $|v_i| \leqslant |v|/2+1$.

By isomorphism the same property is true for a reduced word in the
alphabet $a, \tilde{b}, \tilde{c}, \tilde{d}$ determining an element
in $\tilde{G}$ (and the core in this case is $\widetilde
{\mathcal{N}_0}$ $= \{ 1, \tilde{b}, \tilde{c}, \tilde{d} \}$).

Split the word $w$ as a product of monads $ \ast $ and triads
$x^{-1} \ast x$. If there are at least two monads or at least two
triads we get after rewriting shortening at each of coordinates. The
remaining case is the case of a word of the form $ x^{- 1} \ast x^{}
\ast$ and $\ast x^{- 1} \ast x $ for which one checks that reduction
of length occurs in the second step.

This completes the proof of Proposition~\ref{2.12}.
\end{proof}


>From this proposition we get an algorithm to solve the word problem:
the branch algorithm for $L$. Let us describe it further.

Let $w$ be a word in the letters $a$, $b$, $c$, $d$, $x$. The
problem is to check if $w =1$ in $L$. The notation $w\equiv_L w'$
means that the two elements of $L$ defined by the words $w$ and $w'$
are equal.

1) Reduce $w$ in $\Gamma$. If $w$ is the empty word then in $L$, $w
\equiv_L1$. If it is not the empty word, compute the exponent
$\exp_a w$ (that is the sum of exponents of $a$ in $w$). Check if
this number is even. If NO then $w \not\equiv_L 1$. If YES go to 2).

2) Rewrite $w$ as a pair $( w_0, w_1 )$ using the rewriting map
$\varphi = ( \varphi_0, \varphi_1 )$. Apply 1) successively to $w_0,
w_1$ and follow steps 1) and 2) alternatively. Either, at some step
one obtains a word with odd $exp_a$ or (after $n$ steps) one obtains
that all $2^n$ words represent the identity element in $\Gamma$
(observe that the word problem in $\Gamma$ is solvable by using the
normal form for elements).

Note that $w \equiv_L 1$ $\Leftrightarrow$ ($w_0 \equiv_L 1 $ and
$w_1 \equiv_L 1)$. Applying this procedure 3 times yields either the
answer NO (the elemnt is not the identity) or a set of $8$ words
$w_{i, j, k}$ with $i, j, k \in \{ 0, 1 \}$ which - by
Proposition~\ref{2.12}- are strictly shorter than $w$. This
algorithm solves the word problem.

\begin{lemma}\label{2.14}
Let $w$ be a word in the generators. Let $w \equiv_L ( w_0, w_1 ) \alpha$, $\alpha = a$ or $\alpha = 1$ depending on the parity of
the exponent $\exp_a w$, and the triple $( w_0, w_1 ), \alpha$ is obtained
from $w$ by applying once
 the rewriting process described above. Then
\[\exp_x (w) = \exp_x (w_0) + \exp_x (w_1).\]
\end{lemma}

\begin{proof}
The rewriting process uses the relations $b = (a, c), c = ( a, d ),
d = (1, b)$ and $x = (1, x) a$, $x^{- 1} = ( x^{- 1}, 1 ) a$ which
 do not change the total exponent of $x$. The reduction in group $\Gamma$
 also doesn't change the exponent.
\end{proof}

\begin{lemma}\label{2.15}
The abelianization $L / [ L, L ]$ is infinite. The image of
$x$ in $L / [ L, L ]$ is of infinite order.
\end{lemma}

\begin{proof}
Any element in the commutator group can be expressed as a product of
commutators $[ u, v ]$. Choosing the words in $a, b, c, d, x$
representing $u$ and $v$, we get that any element in $[ L, L ]$ can
be written as a word $w$ with $\exp_x w=0$. Suppose that for some $n
\geqslant 1,$ $x^n \in [ L, L ]$. We get a word $w = x^n \Pi [ u_i,
v_i ]$ in the letters $a, b, c, d, x$ with total exponent $n$ for
$x$ which represents the identity element in $L.$ Choose $w$ of
minimal length with this property. Applying the rewriting process at
most $3$ times to $w$, we get a set of 8 words $w_{i j k}, i,j,k
\in{0,1}$ representing the identity element in $L$ with the sum of
exponents of the symbol $x$ different from zero. Hence at least one
of them has non zero $\exp_x$. The words $w_{i j k}$ are shorter
than $w$, a contradiction.
\end{proof}

The proof of Lemma~\ref{2.15} completes the proof of
Theorem~\ref{2.2}.

\def\cprime{$'$}
\providecommand{\bysame}{\leavevmode\hbox
to3em{\hrulefill}\thinspace}
\providecommand{\MR}{\relax\ifhmode\unskip\space\fi MR }
\providecommand{\MRhref}[2]{%
  \href{http://www.ams.org/mathscinet-getitem?mr=#1}{#2}
} \providecommand{\href}[2]{#2}


\end{document}